\documentstyle[12pt]{amsart}
\theoremstyle{plain}
\newtheorem{Thm}{ }[section]

\title{Subbundles of maximal degree}
\author{Montserrat Teixidor i Bigas}
\address{Mathematics Department, Tufts University, Medford MA
02155, USA} \email{montserrat.teixidoribigas@@tufts.edu}
\begin{document}
\maketitle
\begin{section}{Introduction}
Let $C$ be a curve of genus $g$ and $E$ a generic (semistable)
vector bundle of rank $r$ and degree $d$. Fix a rank $r'<r$ and a
degree $d'$ for subsheaves $E'$ of $E$. If
$r'd-rd'=r'(r-r')(g-1)$, the number of such subbundles is finite.
We shall denote this number with $m(r,d,r',g)$.

The number $m(2,d,1,g)$ had been known for a while (cf.
\cite{S,L,G}. In the case $r'=1$, the number was computed recently
by Oxbury \cite{O} and Okonek-Teleman \cite{OT}. A method for
computing these numbers when $(r',d')$ are relatively prime has
been presented by Lange-Newstead in \cite{LN}. A different method
without any restrictions on $r,d,r'$ is given in \cite{Ho}.

The purpose of this paper is to introduce another approach to the
problem. In the present form it works only for $r'=1$ and for
$r'=2, r=4$. The advantage of the method is that it requires very
little technical background and gives very explicit results.

We need to prove  the following

\begin{Thm}
\label{Thm1} {\bf Theorem} (Oxbury)
 Let $C$ be a curve of genus
$g$. Let $E$ be a generic vector bundle on $C$ of rank $r$ and
degree $d$. Choose $d'$ so that $d-rd'=(r-1)(g-1)$. Then, the
number of line subbundles  of degree $d'$ of $E$ is $r^g$.
\end{Thm}

\begin{Thm}
\label{Thm2} {\bf Theorem}
 Let $C$ be a curve of genus
$g$. Let $E$ be a generic vector bundle on $C$ of rank $4$ and
degree $d$. Choose $d'$ so that $2d-4d'=4(g-1)$. Denote by
$a_g$(resp $b_g$) the number of  subbundles  of rank two and
degree $d'$ of $E$ for $d'$ even (resp odd). Then

$$a_g={g\choose 0}6^{g}+{g\choose 2}6^{g-2}2^2+...+{g\choose
g-\epsilon}6^{\epsilon}2^{g-\epsilon}$$
$$b_g={g\choose 1}6^{g-1}2+{g\choose 3}6^{g-3}2^3+...+{g\choose
g-1+\epsilon}6^{1-\epsilon}2^{g-1+\epsilon}$$

Here $g \equiv \epsilon (2),\ \epsilon \in \{ 0,1\}$

\end{Thm}
\bigskip
The author is a member of the research group "Vector Bundles on
Algebraic Curves".

\end{section}
\begin{section}{Proof of the results}

The proof of the two theorems above will be done by induction on
the genus. We assume the result for a curve of genus $g$. We then
consider a curve of genus $g+1$ obtained by choosing a  curve
$C_1$ of genus one and a  point $P$ on $C_1$ and a  curve $C_g$ of
genus $g$  and a point $Q$ on $C_g$. Glue then $P$ and $Q$ to
obtain a curve $C_{g+1}$ of genus $g+1$. Take a generic vector
bundle $E_{g+1}$ on $C_{g+1}$ and count the number of subbundles
of maximal rank on $C_{g+1}$. We then check that each of these
subbundles corresponds to a non-singular point of the quotient
scheme of $E$ of suitable rank and degree, hence it should be
counted with multiplicity one. This will complete the proof.

{\it Proof of \ref{Thm1}}  We first check the result when $g=1$.
In our situation ($r'=1, g=1$), the numerical condition for the
existence of a finite number of subbundles is $d=rd'$. Then, the
generic vector bundle $E_1$ of rank $r$ and degree $rd'$ on an
elliptic curve is the direct sum of $r$ generic ( and therefore
different) line subbundles of degree $d'$ on $C$ (see \cite{T} ,
Step 3 p.347). A line subbundle of degree $d'$ of  $E_1$ must be
one of the $r$ that appears in the direct sum decomposition. There
are $r$ of them and this agrees with the formula above.

Assume now the result for $g$ and check it for $g+1$ using the
reducible curve described above. Choose numbers $r,d,d'$ such that
$$(*)\ d-rd'=(r-1)((g+1)-1).$$
 Take then  a generic vector
bundle $E_1$ of rank $r$ and degree $r-1$ on $C_1$. Take a generic
vector bundle $E_g$ on $C_g$ of rank $r$ and degree $d-r+1$ and
glue them by a generic gluing. This gives a generic vector bundle
on $C_{g+1}$ (again by \cite{T}).

From (*), $(d-r+1)-rd'=(r-1)(g-1)$. Hence (by the genericity of
$E_g$), the largest degree of a line subbundle of $E_g$ is $d'$.
  By the semistability of $E_1$, the largest possible degree of a line
   subbundle of $E_1$ is zero.
Hence, the only way to obtain a line subbundle of $E_{g+1}$ of
degree $d'$ is  by gluing a line subbundle of degree zero of $E_1$
with a line subbundle of degree $d'$ of $E_g$. By induction
assumption, there are $r^g$ line subbundles of degree $d'$ of
$E_g$. We need to compute how many of the line subbundles of
degree zero of $E_1$ glue with one given direction $V_1$ in the
fiber $(E_1)_P$ of $E_1$ at $P$. Consider the exact sequence
$$0\rightarrow E'_1\rightarrow E_1\rightarrow (E_1)_P/V_1\rightarrow
0$$

A subbundle of $E_1$ that glues with the fixed direction $V_1$
gives rise to a subbundle of $E'_1$. As $deg(E'_1)=deg
E_1-(r-1)=0$, there are $r$ such line subbundles. Hence, $E_{g+1}$
has $r^g\times r=r^{g+1}$ subbundles, as claimed.

We need to check now that each of them needs to be counted with
multiplicity one. Equivalently, we need to show that for each such
sublinebundle $L_{g+1}$ of $E_{g+1}$ the quotient
$E_{g+1}/L_{g+1}$ is a non-singular point of the  scheme of
quotients of rank $r-1$ and degree $d-d'$ of $E$. As the set of
such $L_{g+1}$ is finite, the dimension of the quotient scheme is
zero. The tangent space to the quotient scheme at the point
$E_{g+1}/L_{g+1}$ is given by $Hom(L_{g+1},
E_{g+1}/L_{g+1})=H^0(L_{g+1}^*\otimes E_{g+1}/L_{g+1})$. We need
to check that this vector space is zero-dimensional.

From \cite{RT} Claim p.495, the pair $L_g, E_g/L_g$ is generic in
the product of the moduli spaces of vector bundles of rank one and
$r-1$. Then, from Hirschowitz's Theorem (\cite{Hi} 4.6 or
\cite{RT} 1.2), $h^0(L_g^*\otimes E_g/L_g)=0$. Write
$F_1=E_1/L_1$. Then, a section of $Hom(L_{g+1}, E_{g+1}/L_{g+1})$
is a section of $Hom(L_1, F_1))$ that vanishes at $P$.
Equivalently, this is a section of $L_1^*\otimes F_1(-P)$. We need
to show that $h^0(F_1\otimes  L_1^*(-P))=0$. This could be seen
using the results in \cite{Lange}, Lemma 2.5. We provide instead
an ad hoc proof.
 From the genericity
of all the data, it suffices to show that $h^0(F_1\otimes
L_1^*(-P))=0$ for at least one choice of data. Take as $F_1$ a
direct sum of $r-1$ generic line bundles of degree one. Take as
$E_1$ a generic extension
$$0\rightarrow L_1\rightarrow E_1\rightarrow F_1\rightarrow 0$$
We claim that $E_1$ is indecomposable. This is equivalent to
showing that it is semistable. If this were not the case, there
would be a subsheaf $E'$ of $E_1$ contradicting semistability. We
can assume that $E'$ is semistable, otherwise it sufices to
replace it with a direct summand of maximum slope. Then,
$${d_{E'}\over r_{E'}}>{d_E\over r_E}={r-1\over r}=1-{1\over r}$$
Hence, $d_{E'}>r_{E'}-r_{E'}/r$. As $r_{E'}<r$, this implies
$d_{E'}\ge r_{E'}$. Consider then the exact diagram

$$\begin{matrix} 0&\rightarrow &L_1&\rightarrow& E_1&\rightarrow &F_1&\rightarrow &0\\
 & &\uparrow & &\uparrow & & \uparrow&  & \\
 0&\rightarrow &L'&\rightarrow& E'&\rightarrow &F'&\rightarrow &0\\
 \end{matrix}$$

As $L'$ is a subsheaf of $L_1$ and $F'$ a subsheaf of $F_1$ and
both $L_1$ and $F_1$ are semistable, the condition $\mu(E')\ge 1$
implies $L'=0$ and  $E'=F'$ is a direct summand of $F_1$. But this
contradicts the genericity of the extension defining $E_1$.

Now $F_1\otimes (L_1)^*$ is again a direct sum of generic line
bundles of degree one. Then a generic choice of $P$ gives
$h^0(F_1\otimes  L_1^*(-P))=0$ as required. Note that the whole
picture fits together when we take as $V$ any subspace of $E_P$
containing $(L_1)_P$

 This concludes the
proof  of \ref{Thm1}.
\bigskip

{\it  Proof of \ref{Thm2}}. When $g=1$, if $2d-4d'=4(g-1)=0,\
d=2d'$. For odd $d'$, a generic vector bundle of rank four and
degree $2d'$ is the direct sum of two indecomposable vector
bundles of rank two and degree $d'$ (\cite{T} Step3
p.347). Hence,
$b_1=2$. If $d'$ is even ($d'=2\bar d$), then the generic vector
bundle of rank four and degree $4\bar d$ is the direct sum of four
line bundles of degree $\bar d$. The subbundles of rank four and
degree $2\bar d$ are the direct sum of two of these line bundles.
There are ${4\choose 2}=6$ such choices. Hence, $a_1=6$.

Assume now the result for $g$ and prove it for $g+1$. Choose
$d,d'$ such that

$$(*)2d-4d'=4((g+1)-1)$$
Consider the same type of reducible curve of genus $g+1$ as
before. Take $E_1$ a generic vector bundle of rank four and degree
two on $C_1$. Take $E_g$ a generic vector bundle of rank four and
degree $d-2$ on $C_g$.

By the semistability of $E_1$, the maximum degree of a subbundle
of rank two of $E_1$ is one. From (*), $2(d-2)-4d'=4(g-1)$. Hence,
the maximum degree of a subbundle of rank two of $E_g$ is $d'$.
So, a subbundle of rank two and degree $d'$ of $E_{g+1}$ has its
degree split as either $0,d'$ or $1,d'-1$. We point out that there
are no subbundles with degree split as $1,d'$. This follows from
the fact that $E_1$ has only a finite number of subbundles of
degree one and $E_g$ has only a finite number of subbundles of
degree $d'$. As the gluing is generic, they cannot glue with each
other.

As $E_1$ has rank four and degree two, it has $b_1=2$ subbundles
of degree one and rank two. We need to know how many of the
subbundles of rank two and degree $d'-1$ of $E_g$ glue with a
given subspace of dimension two $V_2$ of $E_P$. Consider the exact
sequence
$$ 0\rightarrow E'_g \rightarrow E_g \rightarrow
(E_g)_P/V_2\rightarrow 0$$
 We see that we need to consider subbundles of rank two and degree
 one of $E'_g$. From the genericity of the gluing at $P$, the
 spaces $V_2$ give rise to generic vector spaces of $(E_g)_P$.From
 \ref{Prop} below, $E'_g$ is a generic vector bundle.

  Note that $d-2-2=d-4$ and  from (*) $2(d-4)-4(d'-1)=4(g-1)$.
 As $E'_g$ is a generic vector bundle of rank four and
 degree $d-2$, the number of its
 subbundles of rank two and (maximal) degree $d'-1$ is $b_g$ if $d'$ is even
 and $a_g$ if $d'$ is odd. Hence, the contribution of the
 subbundles with splitting type $1,d'-1$ to $a_{g+1}$ (resp $b_{g+1}$) is
 $2b_g$ (resp $2a_g$).

 Look now at the splitting of the degrees as $0,d'$. We need to
 consider subbundles of $E_1$ that glue with a given subspace of
 dimension two
 $V_2$ of $(E_1)_P$. This is equivalent to considering subbundles
 of rank two and degree zero of $E'_1$ with $E'_1$ defined by the
 exact sequence
$$ 0\rightarrow E'_1\rightarrow E_1 \rightarrow
(E_1)_P/V_2\rightarrow 0.$$

From \ref{Prop} below, $E'_1$ is generic. From the genus one case,
there are $a_1=6$ such subbundles. We then obtain
$$ a_{g+1}=6a_g+2b_g, b_{g+1}=6b_g+2a_g$$
Using this and the values of $a_1, b_1$, one can check the
validity of the expression in \ref{Thm2} by induction on $g$.

The proof that these correspond to non-singular points of the
quotient scheme is essentially the same  as before. We only need
to check that the pair consisting of such a subbundle $F$ of
degree $d'-1$ and the quotient $E_g/F$ are generic. As $F$ is a
subbundle of maximal degree of $E'_g$, the pair $(F, E'_g/F)$ is
generic. Then, from the diagram
$$\begin{matrix}
0&\rightarrow&F&\rightarrow&E'_g&\rightarrow&E'_g/F&\rightarrow&0\\
& &\downarrow & &\downarrow & &\downarrow & & \\
0&\rightarrow&F&\rightarrow&E_g&\rightarrow&E_g/F&\rightarrow&0\\
& & & &\downarrow & &\downarrow & & \\
 & & & &(E_g)_P/V_2&\rightarrow&(E_g)_P/V_2& & \\\end{matrix}$$

 Dualizing the last column, we obtain
 $$0\rightarrow (E_g/F)^*\rightarrow (E'_g/F)^*  \rightarrow W_2 \rightarrow 0$$
 where $W_2$ is a skyscraper sheaf with support on $P$ and fiber of dimension two.
   Then, the genericity of $E_g/F$ follows from \ref{Prop}.

\begin{Thm} \label{Prop}
 {\bf Proposition} Let $C$ be a curve of genus $g$, $E$
a generic vector bundle of rank $r$ and degree $d$. Choose any
point $P$ on $C$ and a generic surjective map $E_P\rightarrow V_k$
where $V_k$ is a $k$-dimensional vector space. Then, the kernel of
the composition of the natural morphism $E\rightarrow E_P$ with
the map above is a generic vector bundle of rank $r$ and degree
$d-k$.
\end{Thm}
{\it Proof} Consider the set $X$ of pairs consisting of a vector
bundle $E$ and a not necessarily surjective map as above. Then,
$X$ is irreducible of dimension $r^2(g-1)+1+rk$. When we require
the map to be surjective, we obtain a non-empty open set in $X$
which is therefore of the same dimension and irreducible. To every
element in $X$, we can associate the kernel of the composition map
$E\rightarrow E_P\rightarrow V_k$. This is a vector bundle $E'$.
Moreover such an $E'$ appears as kernel  in an $rk$ dimensional
family of these vector bundles. In fact, from
$$0\rightarrow E'\rightarrow E\rightarrow V_k\rightarrow 0$$
one gets
$$0\rightarrow E^*\rightarrow (E')^*\rightarrow V_k\rightarrow 0$$
and $E$ and the surjective map can be recovered in this way.
Hence, $E'$ moves in an $r^2(g-1)+1$ dimensional set and is
therefore generic. This concludes the proof of the Proposition.

\end{section}

\end{document}